 \documentclass[12pt]{amsart}
 \usepackage{amsmath,amssymb}
 \usepackage{pstricks}
 \usepackage{pst-node,pst-tree}
 \newtheorem{thm}{Theorem}[section]
 \newtheorem{cor}[thm]{Corollary}
 \newtheorem{lem}[thm]{Lemma}
 \newtheorem{prop}[thm]{Proposition}

 \newtheorem{prob}[thm]{Problem}
 \theoremstyle{definition}
 
 \theoremstyle{remark}
 
 \numberwithin{equation}{section}
 \newcommand{\set}[1]{\left\{#1\right\}}
 \newcommand{\floor}[1]{\left\lfloor#1\right\rfloor}
 \newcommand{\ceil}[1]{\left\lceil#1\right\rceil}
 
\begin{document}
 \title{Chain Graphs have Unbounded Readability}%
 \author{Martin Charles Golumbic}
 \address{Golumbic: Caesarea Rothschild Institute and
 Department of Computer Science, University of Haifa, Israel}
 \email{golumbic@cs.haifa.ac.il}
 \author{Uri N. Peled}
 \address{Peled: The University of Illinois at Chicago, United States}
 \email{uripeled@uic.edu}
 \thanks{UNP and UR thank the Caesarea Edmond Benjamin de Rothschild Foundation
  Institute for Interdisciplinary Applications of Computer Science
  at the University of Haifa, Israel, for partial support, and
  Daniel Kobler for many discussions concerning
  Problem~\ref{trianglefree}. The authors thank Gyuri Turan for the
  material in Subsection~\protect\ref{subs:read}.}
 \author{Udi Rotics}
 \address{Rotics: Netanya Academic College, Israel}
 \email{rotics@mars.netanya.ac.il}
 \keywords{Read-$k$ graphs}%
 \date{October 11, 2006}%
 \begin{abstract}
   A triangle-free graph $G$ is called read-$k$ when there exists a
   monotone Boolean formula $\phi$ whose variables are the vertices
   of $G$ and whose minterms are precisely the edges of $G$, such
   that no variable occurs more than $k$ times in $\phi$. The
   smallest such $k$ is called the readability of $G$. We exhibit a
   very simple class of bipartite chain graphs on $2n$ vertices with
   readability $\Omega\left(\sqrt{\frac{\log n}{\log \log n}}\right)$.
 \end{abstract}
 \maketitle
 \section{Introduction}
 \subsection{Terminology}\label{subs:term}
  We consider monotone Boolean formulas --- formulas for short ---
  i.e., formulas $\phi$ built from variables $a_1,\ldots, a_n$ using
  the Boolean operations $\vee$ and $\wedge$, which we denote as $+$
  and $*$ for convenience. If no variable appears more than $k$
  times in $\phi$, we say that $\phi$ is \emph{read-$k$}. A monotone
  Boolean function $F$ is said to be \emph{read-$k$} if $F$ has a
  logically equivalent read-$k$ formula. The \emph{readability} of a
  monotone Boolean function $F$ is the smallest $k$ such that $F$ is
  read-$k$. In general determining the readability of a monotone
  Boolean function might be quite difficult, since to the best of
  our knowledge it is not known whether there is a polynomial-time
  algorithm which, given a monotone Boolean function $F$ in an
  irredundant DNF or CNF representation, decides whether or not $F$
  has a read-$k$ formula, for fixed $k \geq 2$.

  Given a formula $\phi$, we can, using distributivity and
  idempotency, write a formula logically equivalent to $\phi$ in the
  form of sum of products of distinct variables, which we call the
  complete sum of products of $\phi$, denoted by
  $\text{CSOP}(\phi)$. Using the absorption rule $\alpha +
  \alpha*\beta \equiv \alpha$ we can simplify $\text{CSOP}(\phi)$ by
  eliminating products containing other products, obtaining the sum
  of minterms of $\phi$, denoted by $\text{SOP}(\phi)$. Each formula
  $\phi'$ logically equivalent to $\phi$ satisfies
  $\text{SOP}(\phi') = \text{SOP}(\phi)$, so we denote it by
  $\text{SOP}(F)$, where $F$ is the Boolean function given by
  $\phi$. For example, $\phi = a_1*(a_1+a_2)$ is read-$2$,
  $\text{CSOP}(\phi) = a_1 + a_1*a_2$, and $\text{SOP}(\phi) = a_1$.

  With every monotone Boolean function $F$ on the variables
  $a_1,\ldots,a_n$ we associate a simple graph $G_F$ on the vertex
  set $\set{a_1,\ldots,a_n}$ whose edges are the unordered pairs
  $a_ia_j$ such that $a_i$ and $a_j$ occur in the same term of
  $SOP(F)$. Thus each term of $SOP(F)$ induces a clique in $G_F$.
  For example for $F_1 = a_1*a_2*a_3$ and $F_2 = a_1*a_2 + a_2*a_3 +
  a_3*a_1$, both $G_{F_{1}}$ and $G_{F_{2}}$ are the triangle on
  $\set{a_1,a_2,a_3}$. In the other direction, with every simple
  graph $G$ we associate a formula $\phi(G)$, which is the SOP
  formula whose terms are the maximal cliques of $G$. Thus if $G$ is
  the triangle on $\set{a_1,a_2,a_3}$, then $\phi(G) = F_1$. A
  monotone Boolean function $F$ is said to be \emph{normal} when
  $SOP(F) = \phi(G_F)$. If $G$ is triangle-free, then $\phi(G)$ is
  automatically normal. In that case we say that $G$ is
  \emph{read-$k$} if $\phi(G)$ is read-$k$, and a read-$k$ formula
  for $\phi(G)$ with the smallest possible $k$ is said to be
  \emph{read-optimal} for $G$. This smallest $k$ is called the
  \emph{readability} of $G$.

  For example, if $G$ is a complete bipartite graph $G$ with edges
  $a_ib_j$, then $\phi(G)$ has the read-$1$ formula $(a_1 + \cdots +
  a_m) * (b_1 + \cdots + b_n)$. It follows that if the edges of a
  triangle-free graph $G$ can be covered by complete bipartite
  subgraphs in such a way that each vertex belongs to at most $k$ of
  them, then $G$ is read-$k$.

  We illustrate these concepts on grid graphs. It is well-known
  (see for example \cite{Gur:77,Gur:91}) that a monotone Boolean
  function $F$ is read-$1$ if and only if $F$ is normal and $G_F$
  is a cograph, i.e., $G_F$ does not have a path on 4 vertices as
  an induced subgraph. Since grid graphs are triangle-free but are
  not cographs (unless the grid is 1 by 1), they are not read-$1$.
  On the other hand, it is easy to cover the edges of a grid graph
  $G$ by complete bipartite subgraphs of the form $K_{2,2}$,
  $K_{1,1}$ and $K_{1,2}$ in such a way that each vertex belongs
  to at most two subgraphs. To do this, color the squares of $G$
  with black and white as in Chess, and for each black square take
  its bounding cycle. These $K_{2,2}$ subgraphs cover all the
  internal edges of G. Then cover the uncovered boundary edges
  with $K_{1,1}$ and $K_{1,2}$. This shows that the readability of $G$ is
  $2$.

  \begin{prob}
    \label{trianglefree} Is it true that a triangle-free graph $G$
    always has a read-optimal formula obtained by covering the edges
    of $G$ with complete bipartite subgraphs?
  \end{prob}

  \subsection{Background on readability}\label{subs:read}
  We are indebted to G. Turan~\cite{Tur:06} for the following
  background information on readability of monotone normal Boolean
  functions. Recall that a monotone quadratic Boolean function
  $F$ is normal if and only if $G_F$ is triangle-free.

  \begin{prop}
   Almost all $n$-variable monotone quadratic Boolean functions have
   readability $\Omega(\frac{n}{\log n})$.
  \end{prop}
  \begin{proof}
   \mbox{}\\[-\baselineskip] 
   \begin{enumerate}
     \item\label{Qn} Let $Q_n$ be the number of $n$-variable
     monotone quadratic Boolean functions. Since every subgraph of a
     complete bipartite graph $K_{n,n}$ is triangle-free, $\log Q_n \geq c_1
     n^2$ for some constant $c_1 > 0$.
     \item\label{Mnk} Every monotone formula is associated with a
     parse tree, with variables at the leaves, and $+$ and $*$
     internal nodes representing the Boolean operations in the
     formula. The size of the formula is defined as the number of
     nodes in the parse tree. Let $M_{n,s}$ be the number of of
     $n$-variable monotone Boolean formulas of size $s$, and we
     estimate it as follows. The parse tree is an ordered tree,
     and there are $\frac{1}{s}\binom{2s-2}{s-1} \leq 2^{2s}$
     ordered trees with $s$ nodes. The tree has at most $s$
     internal nodes and at most $s$ leaves. Therefore there are at
     most $2^s$ ways to assign $*$ or $+$ to the internal nodes,
     and at most $n^s$ ways to assign the $n$ variables to the
     leaves. Multiplying everything together, we deduce that
     $M_{n,s} \leq 2^{3s}n^s$. Therefore $\sum_{j=0}^s M_{n,j}
     \leq \sum_{j=0}^s 2^{3j}n^j \leq 2^{3s+1}n^s$ for $n \geq 2$,
     and therefore $\log \sum_{j=0}^s M_{n,j} \leq c_2 s \log n$
     for some constant $c_2
     > 0$.
     \item\label{Conc} If $s \leq \frac{c_1}{c_2} \frac{n^2}{\log n}
     - \varepsilon$ for some $\varepsilon > 0$, then by~(\ref{Mnk})
     and (\ref{Qn}) we have
     \[
     \begin{split}
     \log \sum_{j=0}^s M_{n,j} \leq c_2 s \log n \leq c_1 n^2 - \varepsilon c_2 \log n
     \\
     \leq  \log Q_n - \varepsilon c_2 \log n,
     \end{split}\]
     or equivalently $\frac{\sum_{j=0}^s M_{n,j}}{Q_n} \leq
     \frac{1}{n^{\varepsilon c_2}} \to 0$. Therefore among all
     $n$-variable monotone quadratic Boolean formulas, the
     proportion of those of size at most $s$ tends to zero. So with
     probability 1 an $n$-variable monotone quadratic Boolean
     formula has size at least $\frac{c_1}{c_2} \frac{n^2}{\log n}$,
     and therefore readability $\Omega(\frac{n}{\log n})$.
   \end{enumerate}
  \end{proof}

  No such functions are known explicitly, but there are explicit
  $n$-variable monotone quadratic Boolean functions with monotone
  formula size $\Omega(n \log n))$ and thus readability $\Omega(\log
  n))$. To explain this, we use the concept of graph entropy defined
  by K\"{o}rner~\cite{Kor:73}.
  We adopt its definition as presented in Newman and
  Wigderson~\cite{New:95}. The entropy of a discrete random variable
  $Z$ is defined as $H(Z) = -\sum_z p(z) \log_2 p(z)$, and the
  mutual information of two random variables $X,Y$ is defined as
  $I(X,Y) = H(X) + H(Y) - H((X,Y))$. Let $A(G)$ be the set of all
  maximal stable sets of a graph $G=(V,E)$. Define $\mathcal{Q}(G)$
  to be the set of all probability distributions $Q_{XY}$ on $V
  \times A(G)$ such that (a) $Q_{XY}(v,I) = 0$ if $v \notin I$, (b)
  the marginal distribution $Q_X$ of $Q_{XY}$ on $V$ is the uniform
  distribution on $V$. Then the \emph{entropy} of $G$ is defined as
  $H(G) = \min \set{I(X,Y)}$, where the minimum is taken over all
  random variables $X$ and $Y$ that are distributed according to the
  marginal distributions $Q_X$ and $Q_Y$ of some distribution
  $Q_{XY} \in \mathcal{Q}(G)$.

  Now we use the following three facts. (1) K\"{o}rner~\cite{Kor:73}
  proved that every $n$ vertex graph $G$ satisfies $H(G) \geq \log_2
  (\frac{n}{\alpha(G)})$, where $\alpha(G)$ is the maximum size of a
  stable set of $G$. (2) Newman and Wigderson~\cite{New:95} proved
  that if $G$ is an $n$-vertex graph, the monotone Boolean formula
  size of $\phi(G)$ is at least $H(G)n$. (3) Using an explicit
  Ramsey construction, Alon~\cite{Alo:95} gave explicit $n$-vertex
  triangle-free graphs $G_n$ with $\alpha(G_n) =
  O(n^{\frac{2}{3}})$. Applying~(1)--(3) to $G_n$, we obtain that
  the monotone Boolean formula size of $\phi(G_n)$ is $\Omega(n \log
  n)$.

  Since an $n$-vertex bipartite graph $G$ satisfies $\alpha(G) \geq
  \frac{n}{2}$, it cannot satisfy $\alpha(G_n) = O(n^{1-\varepsilon})$
  for any $\varepsilon > 0$. Therefore the argument in the preceding
  paragraph cannot use a bipartite graph instead of Alon's $G_n$.

  Jukna~\cite{Juk:04} proved that every $\set{C_3,C_4}$-free
  graph $G=(V,E)$ has monotone Boolean formula size at least $|E|/2$ and hence readability
  $\Omega(|E|/|V|)$. Such graphs include many explicit bipartite
  graphs, and also the point-line incidence graphs of the projective
  planes, for which $|E| \sim |V|^{\frac{3}{2}}$. Thus the
  readability for such graphs can be as high as $\Omega(\sqrt{n})$.

  \subsection{Results}\label{subs:res}

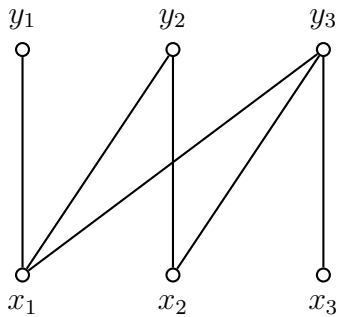
\begin{figure}[ht]
     \begin{center}
      \begin{pspicture}(8,5)
      \psset{radius=.1}
      \Cnode(2,1){x1}
      \nput{270}{x1}{$x_1$}
      \Cnode(4,1){x2}
      \nput{270}{x2}{$x_2$}
      \Cnode(6,1){x3}
      \nput{270}{x3}{$x_3$}
      \Cnode(2,4){y1}
      \nput{90}{y1}{$y_1$}
      \Cnode(4,4){y2}
      \nput{90}{y2}{$y_2$}
      \Cnode(6,4){y3}
      \nput{90}{y3}{$y_3$}
      \ncline{x1}{y1}
      \ncline{x1}{y2}
      \ncline{x1}{y3}
      \ncline{x2}{y2}
      \ncline{x2}{y3}
      \ncline{x3}{y3}
     \end{pspicture}
   \end{center}
  \caption{The graph $G(3)$.}
  \label{Fig:Gn}
  \end{figure}

  The graph $G(n)$ is the bipartite graph with vertices
  $x_1,\ldots,x_n$ and $y_1,\ldots,y_n$ whose edges are the pairs
  $x_iy_j$ with $i \leq j$. Figure~\ref{Fig:Gn} illustrates $G(3)$.

 The graph $G(n)$ is an example of so-called \emph{chain graphs}~\cite{Yan:82},
 also known as \emph{difference graphs}~\cite{Mah:95}. The most general chain
 graph is obtained from $G(n)$ by duplicating vertices, i.e.,
 adding new vertices with the same neighbors as existing vertices.
 It has the same readability as $G(n)$.

 \begin{thm}[\textbf{Main Theorem}]
  \label{mainthm}
  The readability of $G(n)$ is
  \\
  $\Omega\left(\sqrt{\frac{\log n}{\log \log n}}\right)$.
 \end{thm}

 Note that although the lower bound in Theorem~\ref{mainthm} is
 smaller than the ones mentioned above, the graph $G(n)$ is
 bipartite (so is not covered by the arguments of Alon), has $C_4$s
 (so is not covered by the results of Jukna) and has a very simple
 and natural structure. In light of this, Theorem~\ref{mainthm}
 is an interesting result.

 Since $G(n)$ is distance-hereditary, this theorem answers
 affirmatively a question posed in \cite{Gol:95}.

 The following result follows from Theorem~\ref{mainthm}.

 \begin{thm}
  \label{mainthmbigraphs}
  For each $k$, the edges of $G(n)$ cannot be covered by complete
  bipartite subgraphs in such a way that each vertex belongs to at
  most $k$ of them, for sufficiently large $n$.
 \end{thm}

 On the other hand, Theorem~\ref{mainthm} follows from
 Theorem~\ref{mainthmbigraphs} if Problem~\ref{trianglefree} has an
 affirmative answer. We give a graph-theoretical proof of
 Theorem~\ref{mainthmbigraphs} not using Theorem~\ref{mainthm} in
 the Appendix, which may be of independent interest, and served as a
 starting point of our investigations. We also show there that
 $G(n)$ is read-$(1+\ceil{\log_2 n})$.

 Golumbic, Mintz and Rotics~\cite{Gol:95} have shown that if $F$ is
 normal and $G_F$ is a partial $k$-tree, then $F$ is read-$2^k$, and
 thus has bounded readability independent of the number  of vertices
 of $G_F$. Our main theorem continues this line of research with a
 negative result, namely giving a very simple family of bipartite
 graphs with unbounded readability.

 \section{Proof of the Main Theorem}
 We shall be using Greek letters such as $\phi$ and $\psi$ to denote
 formulas. We say that a formula $\psi$ is \emph{as good as} a
 formula $\phi$ when they are logically equivalent and for each
 variable, the number of its occurrences in $\psi$ does not exceed
 the number of its occurrences in $\phi$.

\begin{figure}[ht]
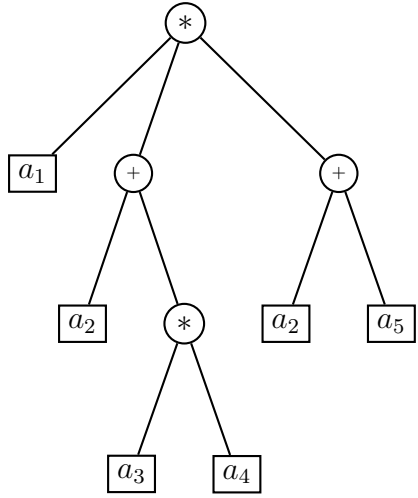

     \begin{center}
      \pstree{\Tcircle{$*$}}
        {
          \Tr{\psframebox{$a_1$}}
          \pstree{\Tcircle{$\scriptscriptstyle +$}}
            {
              \Tr{\psframebox{$a_2$}}
              \pstree{\Tcircle{$*$}}
                {
                  \Tr{\psframebox{$a_3$}}
                  \Tr{\psframebox{$a_4$}}
                }
            }
          \pstree{\Tcircle{$\scriptscriptstyle +$}}
            {
              \Tr{\psframebox{$a_2$}}
              \Tr{\psframebox{$a_5$}}
            }
        }
   \end{center}
  \caption{$\text{tree}(a_1*(a_2+a_3*a_4)*(a_2+a_5))$.}
  \label{Fig:parsetree}
  \end{figure}

 Each formula $\phi$ is associated with a parse tree, denoted by
 $\text{tree}(\phi)$, with the occurrences of the variables of
 $\phi$ at the leaves and the operations $+$ and $*$ of $\phi$ at
 the internal nodes. Figure~\ref{Fig:parsetree} gives an example.

 We can simplify $\text{tree}(\phi)$ by eliminating internal nodes
 corresponding to unary $+$ and $*$ operations, i.e., having a
 single child. Then, using distributivity, we can assume that every
 path down $\text{tree}(\phi)$ alternates between $+$ and $*$ nodes;
 if for example a $+$ node has a $+$ child, remove the child and
 make the grandchildren children of the parent. These operations
 give a logically equivalent formula and do not change the number of
 occurrences of a variable in $\phi$; we always assume they have
 been performed already, as in Figure~\ref{Fig:parsetree}.

 We say that a variable $a_i$ is \emph{isolated} in a formula $\phi$
 if $\phi$ is of the form $\phi = a_i + \psi$.

 A \emph{subformula} of $\phi$ is obtained by taking a node of
 $\text{tree}(\phi)$, removing zero or more of its children but
 leaving at least two children if the node is internal, then taking
 the entire subtree rooted at the resulting node. For example, $a_3$
 and $a_1*(a_2+a_5)$ are subformulas of the formula of
 Figure~\ref{Fig:parsetree}. A subformula $\psi$ of $\phi$ is
 \emph{2-mult} if the root of $\psi$ is a $*$ node and it has
 exactly two children in $\text{tree}(\phi)$. For example, $a_3 *
 a_4$ is a 2-mult subformula of the formula of
 Figure~\ref{Fig:parsetree}, but $a_1*(a_2+a_5)$ is not. A formula
 is said to be \emph{non-redundant} if it does not have a subformula
 of the form $\psi=(a_i + \phi_1)*(a_i + \phi_2)$. Since $a_i +
 \phi_1 * \phi_2$ is as good as $\psi$, every formula $\phi$ can be
 converted to a non-redundant formula that is as good as $\phi$.

 A crucial concept in our proof is that of an extension of $G(n)$. A
 formula $\phi$ is said to be an \emph{extension of $G(n)$} or to
 \emph{extend $G(n)$} when $\text{SOP}(\phi)$ consists of all the
 edges of $G(n)$ (i.e., all the terms of the form $x_i*y_j$ for $1
 \leq i \leq j \leq n$), and in addition zero or more terms, each of
 which is a product of two or more $x_i$ variables or two or more
 $y_j$ variables. For example, $\phi = x_1*(y_1+y_2+y_3) +
 y_3*(x_2+x_3) + x_2*y_2 + x_1 * x_2 * x_3 + y_1 * y_3$ is an
 extension of $G(3)$, but $\psi = x_1*(y_1+y_2+y_3) + y_3*(x_2+x_3)
 + x_2*y_2 + x_2 * y_1 * (x_2 + y_3)$ is not, because
 $\text{SOP}(\psi)$ contains the term $x_2 * y_1$, which is neither
 an edge of $G(3)$ nor a product of two or more $x_i$ or $y_j$
 variables.

 \begin{lem}
  \label{lemma2:claim2} Let $\phi$ be a non-redundant extension of
  $G(m)$. For every edge $x_i*y_j$ of $G(m)$, $\phi$ has a 2-mult
  subformula of the form $(x_i + \phi_1)*(y_j + \phi_2)$.
 \end{lem}

 \begin{proof}
  Since the term $x_i*y_j$ occurs in $\text{SOP}(\phi)$, $\phi$ has
  a subformula of the form $\phi' = (x_i + \phi_1)*(y_j + \phi_2)$
  that contributes this term. If $\phi'$ is 2-mult, we are done. If
  not, this is due to another subformula multiplying $\phi'$ at the
  same level of $\text{tree}(\phi)$, in other words, $\phi$ has a
  subformula of the form $\phi' * \psi$, and because $\phi'$
  contributes $x_i * y_j$ to $\text{SOP}(\phi)$, so does $\phi' *
  \psi$. The formula $\psi$ cannot be a leaf of $\text{tree}(\phi)$,
  because such leaf could only be $x_i$ or $y_j$, and this would
  contradict the non-redundancy of $\phi$. Therefore $\psi$ is
  rooted at a $+$ node or at a $*$ node. In fact we may assume that
  $\psi$ is rooted at a $+$ node, for if $\psi$ has the form
  $\psi=\psi_1 * \psi_2$, we replace $\psi$ with $\psi_1$, and if
  $\psi_1$ still is not rooted at a $*$ node, we continue this
  process of taking the first factor.

  By the non-redundancy of $\phi$, $\psi$ is neither of the
  form $x_i + \psi_1$ nor of the form $y_j + \psi_2$, and therefore
  $\psi$ itself contributes $x_i * y_j$ to $\text{SOP}(\phi)$.

  We now repeat the same argument on $\psi$, and obtain that $\psi$
  has a subformula of the form $\psi' = (x_i + \psi_1)*(y_j +
  \psi_2)$ that contributes the term $x_i*y_j$ to
  $\text{SOP}(\phi)$. If $\psi'$ is 2-mult we are done. If not, we
  notice that because $\phi'$ is rooted at a $*$ node and $\psi$ is
  rooted at a $+$ node, the root of $\psi'$ is a proper descendant
  of the root of $\psi$. Therefore our argument eventually
  terminates in a 2-mult subformula of $\phi$ having the form $(x_i +
  \phi'_1)*(y_j + \phi'_2)$.
 \end{proof}

 We make the notational convention that whenever we write sets of
 the form $\set{i_1,i_2,\ldots,i_n}$ or formulas of the form
 $x_{i(1)} + x_{i(2)} + \cdots + x_{i(n)}$ or $y_{i(1)} + y_{i(2)} +
 \cdots + y_{i(n)}$, we have $i(1) < i(2) < \cdots < i(n)$.

 \begin{lem}
  \label{lemma2:claim3} For every $n$ there exists $m > n$ such that
  every non-redundant read-$k$ extension of $G(m)$ has a subformula
  of the form
  \[(x_{i(1)} + x_{i(2)} + \cdots + x_{i(n)} + \phi_1)*
  (y_{i(1)} + y_{i(2)} + \cdots + y_{i(n)} + \phi_2).\]
  Note that by our notational convention, the subgraph of $G(m)$
  induced by $x_{i(1)},\ldots,x_{i(n)},y_{i(1)},\ldots,y_{i(n)}$ is
  isomorphic to $G(n)$.
 \end{lem}

 \begin{proof}
  Given $n$, we take $m$ as a large enough number, to be specified
  later. Let $\phi$ be a non-redundant read-$k$ extension of $G(m)$.
  By Lemma~\ref{lemma2:claim2}, for each of the edges $x_1*y_j$, $1 \leq j \leq
  m$ of $G(m)$, $\phi$ has a 2-mult subformula of the form
  \[\psi = (x_1 + \phi_1)*(y_j + \phi_2).\]
  We say that $\psi$ \emph{represents the variable $y_j$ with respect to
  $x_1$}. It is possible that a 2-mult subformula $\psi$ of $\phi$
  represents two variables, say $y_{j(1)}$ and $y_{j(2)}$, with respect to
  $x_1$, in which case it has the form
  \[\psi = (x_1 + \phi_1)*(y_{j(1)} + y_{j(2)} + \phi_2).\]
  Since $x_1$ occurs at most $k$ times in $\phi$, there must be at
  least $\ceil{\frac{m}{k}}$ variables $y_{i(1)}, \ldots,
  y_{i\left(\ceil{\frac{m}{k}}\right)}$ among $y_1,\ldots y_m$ all
  represented with respect to $x_1$ by the same 2-mult subformula of
  $\phi$. In other words, $\phi$ has a 2-mult subformula of the form
  \[\psi_1  = (x_1 + \phi_{11})*
  (y_{i(1)} + \cdots + y_{i\left(\ceil{\frac{m}{k}}\right)} + \phi_{12}).\]
  We now consider the variables $x_{i(1)}, \ldots,
  x_{i\left(\ceil{\frac{m}{k}}\right)}$. If at least $n$ of them
  occur isolated in $x_1 + \phi_{11}$, we are done, so we assume
  this is not the case. Therefore at least $n_1 = \ceil{\frac{m}{k}}
  - n$ of these variables (in fact at least $n_1 + 1$ of them), call
  them $x_{j(1)}, \ldots, x_{j(n_1)}$, do not occur isolated in $x_1
  + \phi_{11}$.

  We now repeat the argument for the subgraph of $G(m)$ induced by
  $x_{j(1)}, \ldots, x_{j(n_1)},y_{j(1)}, \ldots, y_{j(n_1)}$.
  Consider the edges $x_{j(1)}*y_{j(l)}$, $1 \leq l \leq n_1$ of
  this subgraph. By Lemma~\ref{lemma2:claim2} and the fact that
  $x_{j(1)}$ occurs at most $k$ times in $\phi$, there is a set of
  $\ceil{\frac{n_1}{k}}$ variables among $y_{j(1)}, \ldots,
  y_{j(n_1)}$, say $y_{i'(1)}, \ldots,
  y_{i'(\ceil{\frac{n_1}{k}})}$, all represented with respect to
  $x_{j(1)}$ by the same 2-mult subformula of $\phi$. In other
  words, $\phi$ has a 2-mult subformula of the form
  \[\psi_2  = (x_{j(1)} + \phi_{21})*
  (y_{i'(1)} +  \cdots + y_{i'\left(\ceil{\frac{n_1}{k}}\right)} + \phi_{22}).\]
  As before, if at least $n$ of the variables $x_{i'(1)}, \ldots,
  x_{i'\left(\ceil{\frac{n_1}{k}}\right)}$ occur isolated in
  $x_{j(1)} + \phi_{21}$, we are done, so we assume this is not the
  case. Therefore at least $n_2 = \ceil{\frac{n_1}{k}} - n$ of these
  variables, call them $x_{j'(1)}, \ldots, x_{j'(n_2)}$, do not
  occur isolated in $x_{j(1)} + \phi_{21}$. And so on.

  If we are not done within $k$ steps, we obtain 2-mult subformulas
  of $\phi$ of the form
  \[\psi_1  = (x_1 + \phi_{11})*
  (y_{i(1)} + \cdots + y_{i(\ceil{\frac{m}{k}})} + \phi_{12}),\]
  with
  \begin{multline*}
   \set{j(1),\ldots,j(n_1)} \subset \set{i(1),\ldots,i\textstyle(\ceil{\frac{m}{k}})} \subset
   \set{1,\ldots,m}, \\
   n_1 = \textstyle\ceil{\frac{m}{k}} - n
  \end{multline*}
  and the variables $x_{j(1)}, \ldots,
  x_{j(n_1)}$ do not occur isolated in
  $x_1 + \phi_{11}$;
  \[\psi_2  = (x_{j(1)} + \phi_{21})*
  (y_{i'(1)} + \cdots + y_{i'(\ceil{\frac{n_1}{k}})} + \phi_{22}),\]
  with
  \begin{multline*}
   \set{j'(1),\ldots,j'(n_2)} \subset \set{i'(1),\ldots,i'\textstyle(\ceil{\frac{n_1}{k}})} \subset
   \set{j(1),\ldots,j(n_1)}, \\
   n_2 = \textstyle\ceil{\frac{n_1}{k}} - n
  \end{multline*}
  and the variables $x_{j'(1)}, \ldots,
  x_{j'(n_2)}$ do not occur isolated in
  $x_{j(1)} + \phi_{21}$;
  \[\psi_3  = (x_{j'(1)} + \phi_{31})*
  (y_{i''(1)} + \cdots + y_{i''(\ceil{\frac{n_2}{k}})} + \phi_{32}),\]
  with
  \begin{multline*}
   \set{j''(1),\ldots,j''(n_3)} \subset \set{i''(1),\ldots,i''\textstyle(\ceil{\frac{n_2}{k}})} \subset
   \set{j'(1),\ldots,j'(n_2)}, \\
   n_3 = \textstyle\ceil{\frac{n_2}{k}} - n
  \end{multline*}
  and the variables $x_{j''(1)}, \ldots,
  x_{j''(n_3)}$ do not occur isolated in
  $x_{j'(1)} + \phi_{31}$;
  And so on. In the general case we use the notation
  $i^{(1)},i^{(2)},\ldots$ for $i',i'', \ldots$ and similarly for
  $j$, and after $k$ steps we obtain
  \[\psi_k  = (x_{j^{(k-2)}(1)} + \phi_{k1})*
  (y_{i^{(k-1)}(1)} + \cdots + y_{i^{(k-1)}(\ceil{\frac{n_{k-1}}{k}})} + \phi_{k2}),\]
  with
  \begin{multline*}
   \set{j^{(k-1)}(1),\ldots,j^{(k-1)}(n_k)} \subset
   \set{i^{(k-1)}(1),\ldots,i^{(k-1)}\textstyle(\ceil{\frac{n_{k-1}}{k}})}
   \\
   \subset \set{j^{(k-2)}(1),\ldots,j^{(k-2)}(n_{k-1})}, \qquad
   n_k = \textstyle\ceil{\frac{n_{k-1}}{k}} - n
  \end{multline*}
  and $x_{j^{(k-1)}(1)}, \ldots, x_{j^{(k-1)}(n_k)}$ do not occur
  isolated in $x_{j^{(k-2)}(1)} + \phi_{k1}$;

  Each of the variables $y_{i^{(k-1)}(1)}, \cdots,
  y_{i^{(k-1)}(\ceil{\frac{n_{k-1}}{k}})}$ occurs in all the
  subformulas $\psi_1,\ldots,\psi_k$. We show that these $k$
  subformulas are distinct, and therefore each of the above
  variables already occurs $k$ times in $\phi$.

  For example, we assume that $\psi_1 = \psi_2$ and obtain a
  contradiction (the argument is the same for $\psi_i = \psi_j$ for
  $i < j$). Let us denote
  \begin{align*}
   \psi_{1L} &= x_1 + \phi_{11}
   \\
   \psi_{1R} &= y_{i(1)} + \cdots + y_{i(\ceil{\frac{m}{k}})} + \phi_{12}
   \\
   \psi_{2L} &= x_{j(1)} + \phi_{21}
   \\
   \psi_{2R} &= y_{i'(1)} + \cdots + y_{i'(\ceil{\frac{n_1}{k}})} + \phi_{22}
  \end{align*}
  Thus $\psi_1 = \psi_{1L} * \psi_{1R}$ and $\psi_2 = \psi_{2L} *
  \psi_{2R}$. By the definition of $\psi_2$, the variable $x_{j(1)}$
  does not occur isolated in $\psi_{1L}$, but it does occur isolated
  in $\psi_{2L}$. Therefore $\psi_{1L} \neq \psi_{2L}$. Since
  $\psi_1$ and $\psi_2$ are 2-mult (they can be factored in only one
  way into two subformulas, up to order), the equality $\psi_{1L} *
  \psi_{1R} = \psi_{2L} * \psi_{2R}$ then implies that $\psi_{1L} =
  \psi_{2R}$ and $\psi_{1R} = \psi_{2L}$. From $\psi_{1L} =
  \psi_{2R}$ it follows that $y_{i'(1)}$ occurs isolated in
  $\psi_{1L}$, and since
  $\set{i'(1),\ldots,i'(\ceil{\frac{n_1}{k}})} \subset
  \set{i(1),\ldots,i(\ceil{\frac{m}{k}}}$, this variable also occurs
  isolated in $\psi_{1R}$. Therefore $\psi_1$ has the form
  $(y_{i'(1)} + \phi_1)*(y_{i'(1)} + \phi_2)$, and this contradicts
  the assumption that $\phi$ is non-redundant. This contradiction
  proves $\psi_1 \neq \psi_2$.

  We have shown that each of the variables
  \[y_{i^{(k-1)}(j)}, \quad 1 \leq j \leq i^{(k-1)}\textstyle(\ceil{\frac{n_{k-1}}{k}})\]
  already occurs $k$ times in $\phi$. We now show that each of the
  variables
  \[x_{i^{(k-1)}(j)}, \quad 1 \leq j \leq i^{(k-1)}\textstyle(\ceil{\frac{n_{k-1}}{k}})\]
  occurs isolated in $x_{j^{(k-2)}(1)} + \phi_{k1}$. We assume that
  for some $1 \leq j \leq i^{(k-1)}(\ceil{\frac{n_{k-1}}{k}})$, the
  variable $x_{i^{(k-1)}(j)}$ does not occur isolated in
  $x_{j^{(k-2)}(1)} + \phi_{k1}$, and obtain a contradiction. By
  construction, this variable also does not appear isolated in any
  of $x_1 + \phi_{11}$, $x_{j(1)} + \phi_{21}$, \ldots,
  $x_{j^{(k-3)}(1)} + \phi_{k-1,1}$. Therefore none of the $k$
  occurrences of the variable $y_{i^{(k-1)}(j)}$ in
  $\psi_1,\psi_2,\ldots,\psi_k$ contributes the term
  $x_{i^{(k-1)}(j)}*y_{i^{(k-1)}(j)}$ to $\text{SOP}(\phi)$. Since
  there are no other occurrences of $y_{i^{(k-1)}(j)}$ in $\phi$,
  the edge $x_{i^{(k-1)}(j)}*y_{i^{(k-1)}(j)}$ of $G(m)$ does not
  occur in $\text{SOP}(\phi)$, contradicting the assumption that
  $\phi$ extends $G(m)$. This contradiction confirms that all of the
  variables
  \[x_{i^{(k-1)}(j)}, \quad 1 \leq j \leq i^{(k-1)}\textstyle(\ceil{\frac{n_{k-1}}{k}})\]
  occur isolated in $x_{j^{(k-2)}(1)} + \phi_{k1}$. We conclude that
  $\psi_k$ is of the form
  \begin{multline*}
  \psi_k = (x_{i^{(k-1)}(1)} + \cdots + x_{i^{(k-1)}\left(\ceil{\frac{n_{k-1}}{k}}\right)}
  + \phi')*
  \\
  (y_{i^{(k-1)}(1)} + \cdots + y_{i^{(k-1)}\left(\ceil{\frac{n_{k-1}}{k}}\right)} +
  \phi_{k2}).
  \end{multline*}
  To conclude the proof, we need only choose $m$ so large that
  $\ceil{\frac{n_{k-1}}{k}} \geq n$.
  We have
  \begin{align*}
    n_1 &\geq \frac{m}{k} - n
    \\
    n_2 &\geq \frac{n_1}{k} - n
    \\
    &\cdots
    \\
    n_{k-1} &\geq \frac{n_{k-2}}{k} - n.
  \end{align*}
  Therefore
  \begin{multline*}
   n_{k-1} \geq \frac{m}{k^{k-1}} - \frac{n}{k^{k-2}} - \cdots - \frac{n}{k} - n
   > \frac{m}{k^{k-1}} - n\left(1+\frac{1}{k}+\frac{1}{k^2} + \cdots\right)
   \\
   = \frac{m}{k^{k-1}} - \frac{nk}{k-1}  \geq \frac{m}{k^{k-1}} - nk.
  \end{multline*}
  It follows that if $m \geq 2nk^k$, we have $\frac{n_{k-1}}{k}> n$, as
  required.
 \end{proof}

 \begin{lem}
  \label{lemma2}
  For every $n$ there exists $m > n$ such that every non-redundant
  read-$k$ extension $\phi$ of $G(m)$ has a subformula of the form
  \[\phi' = (x_{i(1)} + x_{i(2)} + \cdots + x_{i(n)} + \phi_1)*
  (y_{i(1)} + y_{i(2)} + \cdots + y_{i(n)} + \phi_2)\]
  with the following property: Let $\psi$ denote the formula
  obtained from $\phi$ by substituting a new variable $z$ for
  $\phi'$. Then $\text{SOP}(\psi)$ does not contain terms of the
  form $z*x_{i(j)}$ or $z*y_{i(j)}$ for $1 \leq j \leq n$.
 \end{lem}
 \begin{proof}
  We apply Lemma~\ref{lemma2:claim3} for $n+2$ and conclude that
  there exists $m > n+2$ such that every non-redundant read-$k$
  extension of $G(m)$ has a subformula of the form
  \[\phi' = (x_{i(1)} + \cdots + x_{i(n+2)} + \phi_1) * (y_{i(1)} + \cdots + y_{i(n+2)} +
  \phi_2).
  \]
  Define new indices $j(1) = i(2)$, $j(2) = i(3)$, \ldots, $j(n) =
  i(n+1)$, so that $\phi'$ takes the form
  \[\phi' = (x_{j(1)} + \cdots + x_{j(n)} + \phi'_1) * (y_{j(1)} + \cdots + y_{j(n)} +
  \phi'_2),
  \]
  where $\phi'_1 = x_{i(1)} + x_{i(n+2)} + \phi_1$ and $\phi'_2 = y_{i(1)} +
  y_{i(n+2)} + \phi_2$.

  We assume that for some $1 \leq s \leq n$ the term $z*x_{j(s)}$
  occurs in $\text{SOP}(\psi)$ and obtain a contradiction. Replacing $z$
  with $\phi'$ and expanding $\phi'$, we obtain a term $y_{i(1)}*
  x_{j(s)}$ in $\text{CSOP}(\phi)$. This term remains in
  $\text{SOP}(\phi)$, because the latter does not have terms of the
  form $y_{i(1)}$ or $x_{j(s)}$ that could absorb $y_{i(1)}*
  x_{j(s)}$, since $\phi$ is an extension of $G(m)$. Again, since
  $\phi$ is an extension of $G(m)$, we obtain that $y_{i(1)}*
  x_{j(s)}$ is an edge of $G(m)$, a contradiction.

  Similarly no term of the form $z*y_{j(s)}$ occurs in
  $\text{SOP}(\psi)$.
 \end{proof}

 \begin{lem}
  \label{lemma3}
  Suppose $G(n)$ has a read-$k$ extension $\phi$ having a
  subformula of the form
  \[\phi' = (x_1 + x_2 + \cdots + x_n + \phi_1)*(y_1 + y_2 + \cdots + y_n + \phi_2)\]
  with the following property: Let $\phi''$ denote the formula
  obtained from $\phi$ by substituting a new variable $z$ for
  $\phi'$. Then $\text{SOP}(\phi'')$ does not contain terms of the
  form $z*x_i$ or $z*y_j$.

  Then $G(n)$ has a read-$(k-1)$ extension.
 \end{lem}

 \begin{proof}
  We call a minterm that is a product of both
  $x$ and $y$ variables \emph{mixed}. So by definition, the mixed
  minterms of an extension of $G(n)$ are precisely the edges of $G(n)$.

  Let $\psi$ be the formula obtained from $\phi$ by substituting $1$
  (i.e., a true value) for $\phi'$. Since each variable
  $x_1,\ldots,x_n,y_1,\ldots,y_n$ occurs in $\phi'$, each variable
  occurs in $\psi$ less often than in $\phi$. Therefore $\psi$ is
  read-$(k-1)$. To complete the proof, we will show that $\psi$
  extends $G(n)$.

  \emph{Assertion 1:} The term $z$ does not occur in
  $\text{SOP}(\phi'')$, for otherwise we expand $z$ and obtain the
  term $x_2*y_1$ in $\text{CSOP}(\phi)$. This term remains in
  $\text{SOP}(\phi)$ because $\phi$ extends $G(n)$, but this implies
  that $G(n)$ has the edge $x_2*y_1$, a contradiction.

  \emph{Assertion 2:} No terms of the form $x_i$ or $y_j$ occur in
  $\text{SOP}(\psi)$. We assume for example that the term $x_i$ occurs
  in $\text{SOP}(\psi)$ and obtain a contradiction.
  Since $x_i$ is in $\text{SOP}(\psi)$, it follows that the term
  $x_i$ or the term $z*x_i$ is in $\text{SOP}(\phi'')$. The
  hypothesis rules out the latter, so the former holds. But this
  implies that $x_i$ is in $\text{SOP}(\phi)$, which contradicts the
  assumption that $\phi$ extends $G(n)$.

  \emph{Assertion 3:} All the mixed terms of $\text{SOP}(\psi)$ are
  quadratic, i.e., of the form $x_i*y_j$. We suppose that a
  non-quadratic mixed term $A$ occurs in $\text{SOP}(\psi)$ and
  obtain a contradiction. Either $A$ or $z*A$ occurs in
  $\text{SOP}(\phi'')$.

  The first case is that $A$ occurs in $\text{SOP}(\phi'')$. Since
  $\phi$ extends $G(n)$, $A$ does not occur in $\text{SOP}(\phi)$.
  Therefore $A$ is absorbed by a proper subterm $B$ occurring in
  $\text{SOP}(\phi)$. This $B$ does not occur in
  $\text{SOP}(\phi'')$, or else it would also absorb $A$ in
  $\text{SOP}(\phi'')$. It follows that $B$ is obtained in
  $\text{SOP}(\phi)$ by multiplying some term of
  $\text{CSOP}(\phi')$ with some subterm $B'$ of $B$. It follows
  that some subterm of $B'$ occurs in $\text{SOP}(\psi)$. Since $B'$
  is a proper subterm of $A$, $A$ does not appear in
  $\text{SOP}(\psi)$, a contradiction.

  The second case is that $z*A$ occurs in $\text{SOP}(\phi'')$. By
  the forms of $\phi'$ and $A$ we have $\phi'*A=A$. Therefore we see
  that after substituting $\phi'$ for $z$, some subterm $B$ of $A$
  occurs in $\text{SOP}(\phi)$. $B$ must be a proper subterm of $A$
  since $\phi$ extends $G(n)$, and thus all mixed terms of
  $\text{SOP}(\phi)$ are quadratic. Then either $B$ or $zB'$ with
  $B'$ a subterm of $B$ occurs in $\text{SOP}(\phi'')$, and in both
  cases a subterm of $B$ occurs in $\text{SOP}(\psi)$. Since $B$ is
  a proper subterm of $A$, $A$ cannot occur in $\text{SOP}(\psi)$, a
  contradiction.

  \emph{Assertion 4:} $\text{SOP}(\phi)$ and $\text{SOP}(\psi)$ have
  the same mixed terms.

  Let $A$ be a mixed term occurring in $\text{SOP}(\phi)$. Then
  $A$ has the form $x_i*y_j$. The first case is that $A$ occurs in
  $\text{SOP}(\phi'')$. In this case a subterm $B$ of $A$ occurs in
  $\text{SOP}(\psi)$, but $B$ cannot be a proper subterm of $A$ by
  Assertion~2, so $A$ occurs in $\text{SOP}(\psi)$. The second case
  is that $A$ does not occur in $\text{SOP}(\phi'')$. In that case
  $A$ appears in $\text{SOP}(\phi)$ as a result of multiplying
  $\phi'$ by some other formulas. Thus $\text{SOP}(\phi'')$ has a
  term $z*B$ where $B$ is a subterm of $A$. This $B$ cannot be a
  proper subterm of $A$ by Assertion~1 and the hypothesis that
  $z*x_i$ and $z*y_j$ do not occur in $\text{SOP}(\phi'')$.
  Therefore $B=A$ and $z*A$ occurs in $\text{SOP}(\phi'')$.
  Substituting $z=1$ we see that a subterm of $A$ occurs in
  $\text{SOP}(\psi)$, and this subterm must be $A$ itself by
  Assertion~2.

  Conversely, let $A$ be a mixed term occurring in
  $\text{SOP}(\psi)$. By Assertion~3 $A$ must be quadratic, i.e.,
  $A$ has the form $x_i*y_j$. The first case is that $A$ occurs in
  $\text{SOP}(\phi'')$. In this case a subterm of $A$ occurs in
  $\text{SOP}(\phi)$, and this subterm must be $A$ itself because
  $\phi$ extends $G(n)$. The second case is that $A$ does not occur
  in $\text{SOP}(\phi'')$. In that case the term $z*A$ occurs in
  $\text{SOP}(\phi'')$. Substituting $\phi'$ for $z$ we see that the
  terms of $\text{CSOP}(\phi'*A)$ occur in $\text{CSOP}(\phi)$. But
  by the forms of $\phi'$ and $A$ we have $\phi'*A=A$. Therefore a
  subterm of $A$ occurs in $\text{SOP}(\phi)$. Again, by the form of
  $A$ and the hypothesis that $\phi$ extends $G(n)$, this subterm is
  $A$ itself.

  We have proven Assertion~4, and therefore, since $\phi$ extends
  $G(n)$, so does $\psi$, as required.
 \end{proof}

 \begin{thm}
  \label{lemma1.1}
  If $G(n)$ has no read-$(k-1)$ extension, then there exists $m > n$
  such that $G(m)$ has no read-$k$ extension.
 \end{thm}
 \begin{proof}
  Suppose the conclusion of the theorem fails, i.e., for each $m >
  n$, $G(m)$ has a read-$k$ extension. Let $m > n$ be the value
  given by Lemma~\ref{lemma2} for $n$. By our supposition $G(m)$ has
  a read-$k$ extension $\rho$. We can find a non-redundant formula
  $\phi$ that is as good as $\rho$. In particular $\phi$ is
  read-$k$, and $\text{SOP}(\phi) = \text{SOP}(\rho)$, so that
  $\phi$ is also an extension of $G(m)$. By Lemma~\ref{lemma2},
  $\phi$ has a subformula of the form
  \[\phi' = (x_{i(1)} + x_{i(2)} + \cdots + x_{i(n)} + \phi_1)*
  (y_{i(1)} + y_{i(2)} + \cdots + y_{i(n)} + \phi_2)\]
  with the following property: Let $\phi''$ denote the formula
  obtained from $\phi$ by substituting a new variable $z$ for
  $\phi'$. Then $\text{SOP}(\phi'')$ does not contain terms of the
  form $z*x_{i(j)}$ or $z*y_{i(j)}$ for $1 \leq j \leq n$.

  Let $\psi$ denote the formula obtained from $\phi$ by substituting
  zero (i.e., false) for all variables except
  $x_{i(1)},\ldots,x_{i(n)},y_{i(1)},\ldots,y_{i(n)}$
  and renumbering $i(1),\ldots,i(n)$ as $1,\ldots,n$. Then $\psi$
  is read-$k$. Since $\phi$ extends $G(m)$, the mixed terms of
  $\text{SOP}(\phi)$ are precisely the edges of $G(m)$. Only the
  edges induced by $x_1,\ldots,x_n$ and
  $y_1,\ldots,y_n$ (in the new numbering) survive the substitution, and these
  edges form $G(n)$. No new non-mixed terms
  appear as the result of the substitution. Therefore $\psi$ extends
  $G(n)$.

  let $\psi'$ be obtained from $\phi'$ by the same substitution and
  renumbering. Then $\psi'$ is a subformula of $\psi$ of the form
  \[\psi' = (x_1 + x_2 + \cdots + x_n + \psi_1)*(y_1 + y_2 + \cdots + y_n + \psi_2)\]
  with the following property: Let $\psi''$ denote the formula
  obtained from $\psi$ by substituting a new variable $z$ for
  $\psi'$. Then $\text{SOP}(\psi'')$ does not contain terms of the
  form $z*x_j$ or $z*y_j$ for $1 \leq j \leq n$. Indeed, suppose
  $z*x_j$ occurs in $\text{SOP}(\psi'')$. Since it does not occur in
  $\text{SOP}(\phi'')$, a proper subterm, i.e., either $z$ or $x_j$,
  occurs in $\text{SOP}(\phi'')$. It follows that either a subterm
  of $x_2*y_1$ or the term $x_j$ occurs in $\text{SOP}(\phi)$, which
  is impossible since $\phi$ extends $G(m)$.

  We have shown that $\psi$ and $\psi'$ satisfy the hypothesis
  of Lemma~\ref{lemma3}, so by its conclusion $G(n)$ has a
  read-$(k-1)$ extension, contradicting the hypothesis of the
  theorem.
 \end{proof}

 \begin{cor}
  \label{Gnunboundedextension} For each $k$, $G(m)$ has no read-$k$
  extension for $m$ sufficiently large.
 \end{cor}
 \begin{proof}
  By Theorem~\ref{lemma1.1} and the fact that $G(2)$ has no read-$1$
  extension, it follows that there exists an $m$ such that $G(m)$
  has no read-$k$ extension. If $G(m+1)$ had a read-$k$ extension,
  we would obtain from it a read-$k$ extension of $G(m)$ by
  substituting zero for for $x_{m+1}$ and $y_{m+1}$.
 \end{proof}

 \begin{cor}
  \label{Gnunboundedreadability}
   For each $k$, $G(m)$ is not read-$k$ for $m$ sufficiently large.
 \end{cor}
 \begin{proof}
  This follows from Corollary~\ref{Gnunboundedextension}, since every
  formula for $G(m)$ is an extension of $G(m)$.
 \end{proof}

 To prove our main theorem, we analyze the proofs above to find out
 how large they require $m$ to be for a given $k$.

 \begin{proof} (of Theorem~\ref{mainthm})
  It follows from the proofs of Lemma~\ref{lemma2:claim3} through
  Corollary~\ref{Gnunboundedextension} that if $G(n)$ has no
  read-$(k-1)$ extension and $m \geq 2nk^k$, then $G(m)$ has no
  read-$k$ extension. Since $G(2)$ has no read-$1$ extension, it
  follows by induction on $k$ that $G(2^k \cdot 1^1  2^2 \cdots
  (k-1)^{k-1})$ has no read-$k$ extension, and therefore it is not
  read-$k$. Since $2^k \cdot 1^1  2^2 \cdots (k-1)^{k-1} \leq 1^1
  2^2 \cdots k^k$, it follows that if $1^1  2^2 \cdots k^k \leq n$,
  then $G(n)$ is not read-$k$. We use the estimate $\log (1^1  2^2
  \cdots k^k) \leq k^2 \log k$. If we substitute $k =
  \floor{\sqrt{\frac{\log n}{\log \log n}}}$, we obtain $k^2 \log k
  \leq \log n$. Therefore for this $k$, $G(n)$ is not read-$k$; in
  other words, the readability of $G(n)$ is
  $\Omega\left(\sqrt{\frac{\log n}{\log \log n}}\right)$.
 \end{proof}

 \section{Appendix}

 We denote by $r_n$ the smallest $k$ such that the edges of $G(n)$
 can be covered by complete bipartite subgraphs in such a way that
 no vertex belongs to more than $k$ subgraphs. Equivalently, $r_n$
 is the smallest number $k$ such that we can give to each vertex of
 $G(n)$ at most $k$ colors in such a way that $x_i$ and $y_j$ share
 a color if and only if $i \leq j$, i.e., if and only if $x_i*y_j$ is
 an edge of $G(n)$. In that case we say that we have
 \emph{represented} $G(n)$ with these colors. The total number of
 colors used does not matter, only how many colors each vertex
 receives. As we mentioned in the Introduction, $r_n$ is an upper
 bound for the readability of $G(n)$.

 \begin{prop}
  \label{nondecreasing}
  $r_n \leq r_{n+1}$.
 \end{prop}
 \begin{proof}
  This follows trivially from the fact that $G(n)$ is an induced
  subgraph of $G(n+1)$.
 \end{proof}

 \begin{lem}
 \label{upperbound}
  $r_{n+m} \leq 1 + r_{\max(n,m)}$.
 \end{lem}
 \begin{proof}
  Assume without loss of generality that $n \leq m$. Consider
  $G(n+m)$. The subgraph $G_1$ induced by $x_1,\ldots,x_n$ and
  $y_1,\ldots,y_n$ is $G(n)$, and the subgraph $G_2$ induced by
  $x_{n+1},\ldots,x_{n+m}$ and $y_{n+1},\ldots,y_{n+m}$ is
  isomorphic to $G(m)$. Let $k = r_m$. We represent $G_2$ with a set
  of colors so that each vertex of $G_2$ receives at most $k$
  colors. Since $r_n \leq k$ by Proposition~\ref{nondecreasing}, we
  can represent $G_1$ by a set of \emph{new} colors so that each
  vertex of $G_1$ receives at most $k$ colors. Since no color is
  common to $G_1$ and $G_2$, we have not represented the
  non-existing edges between $y_1,\ldots,y_n$ and
  $x_{n+1},\ldots,x_{n+m}$. Finally we give a new color to the
  vertices $x_1,\ldots,x_n$ and $y_{n+1},\ldots,y_{n+m}$ to
  represent the edges between $x_1,\ldots,x_n$ and
  $y_{n+1},\ldots,y_{n+m}$. This coloring represents $G(n+m)$ and
  gives at most $k+1$ colors to each vertex.
 \end{proof}

 \begin{cor}
  $r_{2^q} \leq q+1$, or equivalently by
  Proposition~\ref{nondecreasing}, $r_n \leq 1 + \ceil{\log_2 n}$.
 \end{cor}
 \begin{proof}
  This follows from Lemma~\ref{upperbound} and $r_1=1$.
 \end{proof}

 \begin{lem}
   \label{claim4}
   If $r_n \geq k$, then $r_{(2k+1)n} \geq k+1$.
 \end{lem}
 \begin{proof}
   We assume that $r_n \geq k$ but $r_{(2k+1)n} \leq k$ and obtain a
   contradiction. By Proposition~\ref{nondecreasing} we have $k \leq
   r_n \leq r_{(2k+1)n} \leq k$, and consequently
   \[r_n = r_{(2k+1)n} = k.\]
   Let $G = G((2k+1)n)$, and consider a coloring representing $G$
   with at most $k$ colors present at each vertex. We divide up
   $G$ into $2k+1$ induced subgraphs $G_1,G_2,\ldots,G_{2k+1}$
   isomorphic to $G(n)$, $G_i$ being induced by the vertices
   $x_{(i-1)n+1},\ldots,x_{in}$ and $y_{(i-1)n+1},\ldots,y_{in}$,
   $1 \leq i \leq 2k+1$. We call $\set{x_{(i-1)n+1},\ldots,x_{in}}$
   and $\set{y_{(i-1)n+1},\ldots,y_{in}}$ the \emph{opposite sides of $G_i$}.

   The coloring of $G$ also represents $G_i$. This coloring still
   represents $G_i$ if at each vertex of $G_i$ we keep only the
   colors that appear in the opposite side of $G_i$. If the
   resulting coloring has fewer than $k$ colors present at each
   vertex of $G_i$, then $r_n < k$, a contradiction. Therefore $G_i$
   has a vertex with $k$ colors, all appearing in the opposite side
   of $G_i$. We call such a vertex a \emph{distinguished vertex of
   $G_i$}.

   \emph{Assertion 1:} It is impossible that $G_i$ has a
   distinguished vertex $x_p$ and $G_{i+1}$ has a distinguished
   vertex $y_q$. We suppose such distinguished vertices exist and
   obtain a contradiction. The edge $x_p*y_q$ of $G$ necessitates a
   common color to $x_p$ and $y_q$. Since $x_p$ is distinguished,
   this color is present at some vertex $y_r$ of $G_i$, and since
   $y_q$ is distinguished, this color is present at some vertex
   $x_s$ of $G_{i+1}$. This contradicts the non-existence of the
   edge $x_s*y_r$, proving Assertion~1.

   \emph{Assertion 2:} It is impossible that
   $G_i,G_{i+1},\ldots,G_{i+k}$ all have distinguished vertices on
   the same side. Assume for example that $G_j$ has a distinguished
   vertex $y_{d(j)}$ for each $i \leq j \leq i+k$ (the argument is
   similar if $G_i,G_{i+1},\ldots,G_{i+k}$ all have distinguished
   vertices on the $x$ side). Since $y_{d(j)}$ is distinguished, all
   the $k$ colors present at $y_{d(j)}$ appear on the $x$ side of
   $G_j$. Therefore they cannot be present at $y_{d(l)}$ for any $i
   \leq l \leq j-1$, or else a non-existing edge of $G$ would
   appear. It follows that each distinguished vertex $y_{d(j)}$ has
   $k$ colors that are not present at any other distinguished vertex
   $y_{d(j')}$, $j' \neq j$. Now consider the vertex $x_{d(i)}$.
   Since it is adjacent to the $k$ distinguished vertices
   $y_{d(i+1)}, \ldots, y_{d(i+k)}$, it has a common color with each
   of them. This already gives to $x_{d(i)}$ $k$ distinct colors
   that are not present at the distinguished vertex $y_{d(i)}$.
   Since $x_{d(i)}$ has no other colors, the edge $x_{d(i)}*y_{d(i)}$
   is missing, a contradiction. This proves Assertion~2.

   As a consequence of Assertion~1, there exists an index $0 \leq L
   \leq 2k+1$ such that $G_1,\ldots,G_L$ have distinguished vertices
   only on the $y$ side and not on the $x$ side, whereas
   $G_{L+1},\ldots,G_{2k+1}$ have distinguished vertices only on the
   $x$ side and not on the $y$ side. As a consequence of Assertion~2
   we have both $L \leq k$ and $2k+1-L \leq k$, a contradiction,
   which proves the lemma.
 \end{proof}

 Since $r_1 = 1$, Lemma~\ref{claim4} gives $r_{3 \cdot 1} \geq 2$,
 $r_{5 \cdot 3 \cdot 1} \geq 3$, and in general $r_{(2k-1)!!} \geq
 k$, where $(2k-1)!! = (2k-1) \cdot (2k-3) \cdots 3 \cdot 1$. This
 proves Theorem~\ref{mainthmbigraphs}.

 \end{document}